\documentclass[12pt,a4paper,reqno]{amsart}

\usepackage{amssymb}

\newtheorem{thm}{Theorem}
\newtheorem{theorem}{Theorem}
\newtheorem*{theorem*}{Theorem}
\newtheorem{lemma}[thm]{Lemma}

\newtheorem*{defn*} {Definition}

\newtheorem{fremdersatz}{Theorem}

\newenvironment{satOrig}[1]
        {\pagebreak[2] \begin{fremdersatz} {\bf #1} \quad\sl}
        { \end{fremdersatz}}

\newenvironment{proofof}[1]
        {\pagebreak[3] \vspace{12pt}{\bf Proof #1.}  }
        {\hfill $\blacksquare$ \vspace{12pt}}

\def\fa{\mathcal{F}}
\def\nat{{\rm I\! N}}

\def\co{{\mathbb C}}

\def\dk{{\mathbb D}}

\def\real{{\rm Re \,}}
\def\l{\left}
\def\r{\right}
\def\gl{\left\{}
\def\gr{\right\}}
\def\kl{\left(}
\def\kr{\right)}
\def\kj{\overline}
\def\limn{\lim_{n\to\infty}}

\def\In{\subseteq}
\def\mi{\setminus}
\def\abb{\longrightarrow}
\def\plog{\log^+}

\def\beq{\begin{equation}}
\def\eeq{\end{equation}}
\def\beqar{\begin{eqnarray}}
\def\eeqar{\end{eqnarray}}
\def\beqaro{\begin{eqnarray*}}
\def\eeqaro{\end{eqnarray*}}
\def\bsat{\begin{theorem}}
\def\esat{\end{theorem}}
\def\bsator{\begin{satOrig}}
\def\esator{\end{satOrig}}
\def\blem{\begin{lemma}}
\def\elem{\end{lemma}}
\def\bbew{\begin{proofof}}
\def\ebew{\end{proofof}}

\begin{document}

 \title[An Extension of one direction in Marty's normality criterion]
 {An Extension of one direction in Marty's normality criterion}

 \author[J\"urgen Grahl]{J\"urgen Grahl
 }
 \address{J\"urgen Grahl \\
University of W\"urzburg \\
Department of Mathematics \\
97074 W\"urzburg\\
Germany}
\email{grahl@mathematik.uni-wuerzburg.de}

 \author[Shahar Nevo]{Shahar Nevo
 }
 \address{Shahar Nevo\\ Department of Mathematics\\
 Bar-Ilan University, 52900 Ramat-Gan, Israel}
\thanks{$^2$Research of Shahar Nevo was supported by the  Israel Science Foundation,
 Grant No. 395/2007.}
  \email{nevosh@macs.biu.ac.il}


 \begin{abstract}
   We prove the following extension of one direction in Marty's
   theorem: If $k$ is a natural number, $\alpha>1$ and $\fa$ is a
   family of functions meromorphic on a domain $D$ all of whose poles
   have multiplicity at least $\frac{k}{\alpha-1}$, then the normality
   of $\fa$ implies that the family
$$\gl\frac{|f^{(k)}|}{1+|f|^\alpha}\;:\; f\in\fa\gr$$ 
is locally uniformly bounded.  
   \end{abstract}



 \maketitle

 \section{Introduction and main results}\label{introduction}

Our point of departure is the following famous normality criterion of
F. Marty \cite{Marty}. 

\bsator{(Marty's Theorem)} \label{Marty}
A family $\fa$ of meromorphic functions on a
domain $D\In\co$ is normal if and only if the family $\gl f^\# \,:\,
f\in\fa\gr$ of the corresponding spherical derivatives
$f^\#=\frac{|f'|}{1+|f|^2}$ is locally uniformly bounded. 
\esator

In the present paper we investigate the question how normality can be
characterized in terms of the quantity 
$$\frac{|f^{(k)}|}{1+|f|^\alpha} \qquad \mbox{ where } k\in\nat,\;
\alpha>0$$ 
rather than the spherical derivative $f^\#$. A more or less complete
answer is already known for the direction ``$\Longleftarrow$'' in
Marty's theorem (locally uniform boundedness implies normality), but
not for the opposite direction. Hence, we focus our attention on the
latter one.   

But first we summarize the known results concerning direction
``$\Longleftarrow$''.  A substantial (and best possible) improvement
of this direction in Marty's theorem is due to A.~Hinkkanen
\cite{hinkk}: A family of meromorphic (resp.~holomorphic) functions is
already normal if the corresponding spherical derivatives are bounded
on the preimages of a set consisting of five (resp.~three) elements.
(An analogous result for normal functions was earlier proved by P.
Lappan \cite{lappan}.)

As to generalizations of Marty's theorem to higher derivatives,
S.Y.~Li and H.~Xie \cite{LiXie} obtained the following result.

\bsator{} \label{LX} 
Let $k$ be a natural number and $\fa$ a family of functions
meromorphic on a domain $D$ all of whose zeros have multiplicity at
least $k$. Then $\fa$ is normal in $D$ if and only if
$\left\{\dfrac{|f^{(k)}|}{1+|f|^{k+1}}:f\in\fa\right\}$ is
locally uniformly bounded in $D$. The direction ``$\Longrightarrow$''
holds without the assumption on the multiplicities.
\esator

In \cite{NevoPang} a new proof of Theorem \ref{LX} was given which avoids
the use of Nevanlinna theory.

Finally, Y. Xu \cite{Xu} proved the following extension of Hinkkanen's
normality result to higher derivatives. 

\bsator{}\label{Xu}
Let $k$ be a natural number and $\fa$ a family of functions
meromorphic on a domain $D$. Assume that there is a value $w^*\in\co$
and a constant $M<\infty$ such that for each $f\in\fa$ we have
$|f'(z)|+\dots+|f^{(k-1)}(z)|\le M$ whenever $f(z)=w^*$ and that
there exists a set $E\subset\kj\co$ consisting of $k+4$ elements such that
for all $f\in\fa$ and all $z\in D$ we have  
\beq\label{Xu-Cond}
f(z)\in E \quad \implies  \quad\dfrac{|f^{(k)}|}{1+|f|^{k+1}}(z)\le M.
\eeq
Then $\fa$ is a normal family. 
If all functions in $\fa$ are holomorphic, this also holds if one
merely assumes that $E$ has at least $3$ elements. 
\esator

Here and in the following, terms like $\dfrac{|f^{(k)}|}{1+|f|^{k+1}}$
are understood to be continuously extended into the poles of $f$.  Of
course, the use of $\dfrac{|f^{(k)}|}{1+|f|^{k+1}}$ instead of
$|f^{(k)}|$ in Theorem \ref{Xu} is only due to the possibility that
$E$ might contain the point $\infty$; if $\infty\not\in E$, condition
(\ref{Xu-Cond}) can  be replaced by $|f^{(k)}(z)|\le M'$ whenever
$f(z)\in E$ with a suitable constant $M'>M$.     

In ``$\Longleftarrow$'' of Theorem \ref{LX} and in Theorem \ref{Xu},
the condition on the multiplicities of the functions in $\mathcal F$
resp.~the (slightly weaker) condition on the existence of the value
$w^*$ is essential as the non-normal family of polynomials of degree
at most $k-1$ demonstrates. 

So Theorem \ref{Xu} gives a (more or less) complete answer to the
question how direction ``$\Longleftarrow$'' in Marty's theorem can be
extended in terms of $\frac{|f^{(k)}|}{1+|f|^\alpha}$ rather than
$f^\#$. In particular, for arbitrary $\alpha>0$ the locally uniform
boundedness of $\gl \frac{|f^{(k)}|}{1+|f|^\alpha}\;:\;f\in\fa\gr$
implies normality of $\fa$ provided that the zeros of the functions in
$\fa$ appear only with multiplicities at least $k$ -- though the whole 
truth is much stronger since it suffices to investigate the preimages
of ``few'' values. 

 As to the opposite direction in Marty's theorem, we prove the
 following result which generalizes ``$\Longrightarrow$'' in Theorem
 \ref{LX}. 

\begin{theorem}\label{LiXieGeneraliz}
Let $k$ be a natural number, $\alpha>1$ be a real number and let
$\mathcal F$ be a family of functions meromorphic on a domain
$D$ all of whose poles have multiplicity at least
$\frac{k}{\alpha-1}$. Then the normality of $\fa$ implies that 
$${\mathcal  F}_{k,\alpha}:=\left\{\dfrac{|f^{(k)}|}{1+|f|^\alpha}:f\in\mathcal
   F\right\}$$
is locally uniformly bounded.  
\end{theorem}

We explicitly point out two special (and, in some sense, extremal) cases: 
\begin{itemize}
\item[(S1)] 
If $\alpha\ge k+1$ and if $\fa$ is normal, then the conclusion that
$\fa_{k,\alpha}$ is locally uniformly bounded holds without any
further assumptions on the multiplicities of the poles. This is just
the direction ``$\Longrightarrow$'' in Theorem \ref{LX}. (More
precisely, Theorem \ref{LX} settles the case $\alpha=k+1$. But if the
locally uniform boundedness of $\fa_{k,\alpha}$ is proved for 
$\alpha=k+1$, it trivially also holds for $\alpha>k+1$ since $x\mapsto
\frac{1+x^{k+1}}{1+x^\alpha} $ is bounded on $[0,\infty)$ whenever $\alpha>k+1$.)
\item[(S2)] 
If all functions in $\fa$ are holomorphic, then the normality of $\fa$
implies that $\fa_{k,\alpha}$ is locally uniformly bounded for any
$\alpha>1$. For $k=1$, this was already proven in \cite[Theorem~3]{LiuNevoPang}. 
\end{itemize}

In the case $1<\alpha< k+1$ the lower bound $\frac{k}{\alpha-1}$ for the
multiplicities in Theorem \ref{LiXieGeneraliz} is best possible. This
can be seen by considering the single function $f(z)=\frac{1}{z^p}$ with
$p<\frac{k}{\alpha-1}$ near its pole: Here, for a certain $C>0$ we have
$$\dfrac{|f^{(k)}|}{1+|f|^\alpha}(z)\sim C\cdot |z|^{(\alpha-1)
  p-k}\abb\infty \qquad \mbox{ for } z\to 0$$
since $(\alpha-1)p-k<0$. Since $f$ is even zero-free, this example
  also shows that there is no analogue of Theorem \ref{LiXieGeneraliz}
  where the condition on the multiplicities of the poles is replaced
  by a condition on the multiplicities of the zeros.

Even for holomorphic functions the condition $\alpha>1$ cannot be
weakened any further. This is shown by the family of the
functions $f_n(z):=(z-3)^n$ which is normal in the unit disk $\dk$ and
satisfies
\begin{eqnarray*}
\frac{|f_n^{(k)}(z)|}{1+|f_n(z)|^\alpha}
&=&n(n-1)\cdot\ldots\cdot (n-k+1)\cdot \frac{|z-3|^{n-k}}{1+|z-3|^{\alpha n}}\\
&\ge& \frac{1}{2}\cdot(n-k)^k \cdot |z-3|^{n(1-\alpha)-k} \abb\infty  
\qquad (n\to\infty)
\end{eqnarray*}
for all $z\in\dk$ and all $\alpha$ with $0<\alpha\le 1$. 

Theorem \ref{LiXieGeneraliz} is a consequence of the following result
which we hope to be of interest for itself.  To simplify its
statement, we write ``$f_n\overset \chi\Longrightarrow f$ on $D$'' to
indicate that the sequence $\{f_n\}_n$ converges to $f$ uniformly
w.r.t.~the spherical metric on compact subsets of $D$ and
``$f_n\Longrightarrow f$ on $D$'' if the convergence is in the Euclidean
metric.

\begin{theorem}\label{ChenHuaGen}
Let $D$ be a domain in $\co$ and let $k,m,p$ be natural numbers.
\begin{itemize}
\item[(a)] 
Let $\{g_n\}_{n=1}^\infty$ be a sequence of holomorphic functions
$g_n\not\equiv 0$ on $D$ all of whose zeros have multiplicity
at least $m$. If $g_n\Longrightarrow 0$, then 
$$\frac{\kl g_n^{(k)}\kr^m}{g_n^{m-k}}\Longrightarrow 0 .$$
\item[(b)] 
Let $\{f_n\}_{n=1}^\infty$ be a sequence of meromorphic functions
on $D$ all of whose poles have multiplicity at least $p$. If
$f_n\overset\chi\Longrightarrow \infty$, then 
$$\frac{\kl f_n^{(k)}\kr^p}{f_n^{p+k}}\Longrightarrow 0.$$
  \end{itemize}
\end{theorem}

For $p=1$, a proof of (b) was given in \cite{NevoPang}, essentially
based on Weierstra\ss's theorem and induction. For holomorphic
functions (where $p$ can be chosen arbitrarily large) (b) has been
proved by H.~Chen and X.~Hua \cite{chenhua}. Their reasoning (as the
proof of (S2) for $k=1$ in \cite{LiuNevoPang}) uses just Harnack's
inequality (applied to the harmonic functions $\log |f_n|$) and
Cauchy's formula. In the general case considered  here, a much more
careful analysis is required. 

In the case $m\le k$, (a) almost trivially follows from Weierstra\ss's
theorem. So the interesting case in (a) is the case $k<m$. 

We note that the exponents in $\dfrac{\kl g_n^{(k)}\kr^m}{g_n^{m-k}}$
and $\dfrac{\kl f_n^{(k)}\kr^p}{f_n^{p+k}}$ are chosen in such a way
that (under the respective assumptions in (a) and (b)) these functions
are holomorphic and that (for $k<m$) the assumptions on the
multiplicities cannot be weakened without losing the holomorphy.
(Explicit counterexamples are provided by the sequences of the
functions $g_n(z):= z^{m-1}/n$ and $f_n(z):=n/z^{p-1}$ in the unit
disk.) In this sense, Theorem \ref{ChenHuaGen} is best possible.

\section{Proofs}

First let us define some notations. For $z_0\in\mathbb C$ and $r>0,$
we set $\Delta(z_0,r):=\{z\in\co:|z-z_0|<r\} $ and
$\Delta'(z_0,r):=\Delta(z_0,r)\mi\gl z_0\gr. $ Furthermore, we denote
the open unit disk by $\dk:=\Delta(0,1)$. 



The proof of Theorem \ref{ChenHuaGen} is inspired by several ideas used
in the proof of the lemma on the logarithmic derivative (see
\cite{lang}, VI. \S{3}). At some (crucial) point it makes use of the
general form of Poisson-Jensen-Nevanlinna's formula which as a special
case yields the First Fundamental Theorem of Nevanlinna theory. To
simplify notations, we state this application of
Poisson-Jensen-Nevanlinna's formula in terms of a modification of
Nevanlinna theory that was discussed in \cite{grahl-modif}. If $f$ is
a function meromorphic on a disk $\Delta(0,R_0)$ and if $\alpha\in
\Delta(0,R_0)$ is not a pole of $f$, for $|\alpha|<r<R_0$ we define
$$m_\alpha(r,f):=\frac{1}{2\pi}\int_0^{2\pi} \plog |f(r e^{it})|
\cdot \real\frac{re^{it}+\alpha}{re^{it}-\alpha} dt, $$
$$N_\alpha(r,f):=\sum_{|b_k|<r}\log\l|\frac{r^2-\kj{b_k}\alpha}{r(\alpha-b_k)}\r| 
\qquad \mbox{ and } \qquad
T_\alpha(r,f):=m_\alpha(r,f)+N_\alpha(r,f),$$
where the $b_k$ are the poles of $f$, each taken into account
according to its multiplicity. We call $m_\alpha(r,f)$,
$N_\alpha(r,f)$ and $T_\alpha(r,f)$ the {\bf modified proximity
  function, counting function and characteristic} of $f$ with respect
to $\alpha$. Using these quantities, Poisson-Jensen-Nevanlinna's
formula takes the following form \cite[Theorem~1]{grahl-modif}
\beq\label{FFT}
T_\alpha\kl r,\frac{1}{f}\kr=T_\alpha( r,f) + \log\frac{1}{|f(\alpha)|}
\qquad\mbox{ for } |\alpha|<r<R_0
\eeq
provided that $\alpha$ is not a zero or pole of $f$.

Let $n(r,f)$ denote the number of poles of $f$ and $n(r,c,f)$ for
$c\in\co$ the number of poles of $\frac{1}{f-c}$ in the closed disk
$\kj{\Delta(0,r)}$, counted according to their multiplicities. Then
for $|\alpha|<r<R<R_0$ we have the estimate \cite[Lemma~3]{grahl-modif}
\beq\label{diffabsch}
n(r,f)\cdot \frac{(R-r)(R-|\alpha|)}{R^2+r|\alpha|} 
\le N_\alpha(R,f)-N_\alpha(r,f)
\eeq
which is also required in the proof of Theorem \ref{ChenHuaGen}.

\bbew{of Theorem \ref{ChenHuaGen}}
Without loss of generality we may assume that $D=\dk$ is the unit disk
and that the convergence of $\gl g_n\gr_n$ and $\gl f_n\gr_n$  is uniform in $\dk$. 

{\bf I.} First we show that 
$$g_n^k\cdot \l[\kl\frac{g_n'}{g_n}\kr^{(k-1)}\r]^m\Longrightarrow
0.$$
For this purpose we fix $r<R<1$ and set $s:=\frac{1}{2}(r+R)$. 
There exists some $x_0\in(0, 1/e]$ such that for all $x\in
(0,x_0]$ the function $y\mapsto H(x,y)$ where 
$$H(x,y):= \kl\frac{x}{y}\kr^k \cdot \kl m+\frac{1}{(s-r)^2}\cdot \log\frac{y}{x}\kr^{2m}$$
is monotonically decreasing on the interval $[1,\infty)$. 

We consider some fixed function $g\not\equiv0$ holomorphic on $\dk$ all of
whose zeros have multiplicity at least $m$ and which satisfies
$|g(z)|\le x_0$ for all $z\in\dk$.

We define 
$$G_a(z):=\frac{s^2-\kj{a}z}{s(z-a)} \qquad
\mbox{ and } \quad
B:=\prod_{|a_j|<s} G_{a_j}^{m_j}$$
where the $a_j$ are the distinct zeros of $g$ and $m_j\;(\,\ge m)$ their
respective multiplicities. Then 
\beq\label{h-Def}
h:=g\cdot B
\eeq
is holomorphic on $\dk$ and non-vanishing in
$\Delta(0,s)$, and we have 
\beq\label{g-h-B}
g^k\cdot \l[\kl\frac{g'}{g}\kr^{(k-1)}\r]^m
=\frac{h^k}{B^k}\cdot
\l[\kl\frac{h'}{h}\kr^{(k-1)}-\kl\frac{B'}{B}\kr^{(k-1)}\r]^m.
\eeq
From Poisson's formula one easily gets (cf. \cite[Satz 9.2]{jankvolk}) for $|z|<s$
$$\frac{h'}{h}(z) = \frac{1}{2\pi} \int_0^{2\pi} \log |h(se^{it})| 
\cdot\frac{2se^{it}}{(se^{it}-z)^2} \,dt$$
and
$$\l(\frac{h'}{h} \r)^{(k-1)}(z) = \frac{k!}{2\pi} \int_0^{2\pi}\log|h(se^{it})|
\cdot\frac{2se^{it}}{(se^{it}-z)^{k+1}}\,dt.$$
Here, in view of $|B (\zeta)|=1$ for $|\zeta|=s$ we have 
$$|h(\zeta)|=|g(\zeta)|\le x_0 \qquad \mbox{ for } |\zeta|=s, $$
hence by the maximum principle
$$|h(z)|\le x_0\qquad \mbox{ for  } |z|\le s.$$
In particular, $\log |h(\zeta)|<0 $ for $|\zeta|=s$. Therefore, we obtain for $|z|\le r$ 
\beqar\label{h-LogDeriv}
\l|\l(\frac{h'}{h} \r)^{(k-1)}(z)\r| 
&\le& -\frac{2\cdot k!}{2\pi(s-r)^{k+1}} \int_0^{2\pi}\log|h(se^{it})|\,dt\nonumber\\
&=& -\frac{2\cdot k!}{(s-r)^{k+1}} \cdot\log|h(0)|.
\eeqar
To estimate the contribution of $\frac{B'}{B}$ to \eqref{g-h-B}, we use
$$\frac{G'_a}{G_a}(z)=\frac{-\kj{a}}{s^2-\kj{a}z}-\frac{1}{z-a}$$
and
$$\kl\frac{G'_a}{G_a}\kr^{(k-1)}(z)
=(k-1)!\cdot\kl\frac{-\kj{a}^k}{\kl s^2-\kj{a}z\kr^k}+\frac{(-1)^k}{(z-a)^k}\kr$$
and note that for $|z|=r$ and $|a_j|<s$ we have
$|s^2-\kj{a_j}z|\ge s\cdot(s-r)$, hence
$$\l|\kl\frac{G'_{a_j}}{G_{a_j}}\kr^{(k-1)}(z)\r| 
\le  (k-1)!\cdot\kl\frac{1}{(s-r)^k}+\frac{1}{|z-a_j|^k}\kr.$$
So we obtain for $|z|\le r$
\beqar\label{B-LogDeriv}
\l|\kl\frac{B'}{B}\kr^{(k-1)}(z)\r| 
&\le& \sum_{|a_j|<s} m_j\cdot \l|\kl\frac{G'_{a_j}}{G_{a_j}}\kr^{(k-1)}(z)\r| \nonumber\\
&\le&(k-1)!\cdot\kl\frac{n(s,0,g)}{(s-r)^k}+\sum_{|a_j|<s}\frac{m_j}{|z-a_j|^k}\kr
\eeqar
Let's assume that $g$ has at least one zero in $\kj{\Delta(0,s)}$,
i.e. that $n(s,0,g)\ge 1$. Then from (\ref{h-LogDeriv}),
(\ref{B-LogDeriv}) and the (trivial) estimate 
$$a+b+c \le 3abc \quad\mbox{ for all } a,b,c\ge 1$$
we obtain for $|z|\le r$
\beqar\label{SumLeqProd}
&&\l|\kl\frac{h'}{h}\kr^{(k-1)}-\kl\frac{B'}{B}\kr^{(k-1)}\r|(z)\nonumber\\
&\le& \frac{2\cdot k!}{(s-r)^{k+1}} \cdot
\kl \log\frac{1}{|h(0)|} +n(s,0,g)+2^k\cdot \sum_{|a_j|<s}\frac{m_j}{|z-a_j|^k}\kr\nonumber\\
&\le& \frac{6\cdot k!\cdot 2^k}{(s-r)^{k+1}} \cdot\log\frac{1}{|h(0)|}
\cdot n(s,0,g)\cdot \sum_{|a_j|<s}\frac{m_j}{|z-a_j|^k}.
\eeqar
(Here we have used that $-\log|h(0)|\ge -\log x_0\ge 1$ and 
$2^k\cdot \sum_{|a_j|<s}\frac{m_j}{|z-a_j|^k}\ge 1$ since $|z-a_j|\le 2$.)

We fix some $z_0\in \kj{\Delta(0,r)}$ such that $g(z_0)\ne 0$. Then
there is some $j_*=j_*(z_0)$ such that
$|z_0-a_{j_*}|=\min_{|a_j|<s}|z_0-a_j|$. We define
$$B_*:=G_{a_{j_*}}^{m_{j_*}-m}\cdot\prod_{|a_j|<s,\; j\ne j_*}
G_{a_j}^{m_j}
\quad\mbox{ and } \quad 
g_*:=\frac{h}{B_*}.$$
Then $g_*$ is holomorphic on $\kj{\Delta(0,s)}$, and we have 
$$B=B_*\cdot G_{a_{j_*}}^m
\qquad\mbox{ and } \qquad
g_*=g\cdot G_{a_{j_*}}^m .$$
Using $|G_{a_{j_*}}(z)|\le 1$ for $|z|\ge s$, $|g(z)|=|g_*(z)|$ for
$|z|=s$ and the maximum principle we deduce 
$$|g_*(z)|
\le \gl\begin{array}{cll}
|g(z)| & \le x_0 & \mbox{ for } s\le |z|\le R, \\[5pt] 
\max_{|\zeta|=s} |g(\zeta)| & \le x_0 & \mbox{ for } |z|<s,
\end{array}\r.$$ 
i.e.~$|g_*(z)|\le x_0$ for $|z|\le R$, and we obtain 
\beqaro
&&\frac{1}{|B(z_0)|^k}\cdot \kl\sum_{|a_j|<s}\frac{m_j}{|z_0-a_j|^k}\kr^m\\
&\le& \frac{1}{|B_*(z_0)|^k}\cdot\kl\frac{s|z_0-a_{j_*}|}{|s^2-\kj{a_{j_*}}z_0|}\kr^{km}
\cdot (n(s,0,g))^m\cdot \frac{1}{|z_0-a_{j_*}|^{km}}\\
&\le& \frac{1}{|B_*(z_0)|^k}\cdot \frac{1}{(s-r)^{km}}\cdot (n(s,0,g))^m.
\eeqaro
We combine this estimate with \eqref{g-h-B} and \eqref{SumLeqProd} and arrive at
\beqar\label{MainEstim1}
&&\l|g^k\cdot \l[\kl\frac{g'}{g}\kr^{(k-1)}\r]^m\r|(z_0)\nonumber\\
&\le & \l|\frac{h}{B}(z_0)\r|^k\cdot 
\kl\frac{6\cdot k!\cdot 2^k}{(s-r)^{k+1}} \cdot\log\frac{1}{|h(0)|}
\cdot n(s,0,g)\cdot \sum_{|a_j|<s}\frac{m_j}{|z-a_j|^k}\kr^m\nonumber\\
&\le & \l|\frac{h}{B_*}(z_0)\r|^k\cdot 
\kl\frac{6\cdot k!\cdot 2^k}{(s-r)^{2k+1}} \cdot\log\frac{1}{|h(0)|}\kr^m
\cdot (n(s,0,g))^{2m}.
\eeqar
Here, applying the estimate \eqref{diffabsch} and the First Fundamental
Theorem (\ref{FFT}) to $g_*$ and observing that $|g_*(z)|\le1$ for $|z|\le R$
implies $T_{z_0}(R,g_*)=0$, we obtain  
\beqaro
n(s,0,g)&=& m+ n(s,0,g_*)\\
&\le& m+\frac{R^2+s|z_0|}{(R-s)(R-|z_0|)}\cdot N_{z_0}\kl R,\frac{1}{g_*}\kr\\
&\le& m+\frac{2R^2}{(R-s)(R-r)}\cdot \kl T_{z_0}\kl
R,g_*\kr+\log\frac{1}{|g_*(z_0)|}\kr\\
&\le& m+\frac{1}{(s-r)^2}\cdot \log\frac{|B_*(z_0)|}{|h(z_0)|}.
\eeqaro
Inserting this into (\ref{MainEstim1}) and keeping in mind
$|B_*(z_0)|\ge 1$, $|h(z_0)|\le x_0$ and the definition of $x_0$ we
conclude that 
\beqar\label{B*Elim}
&&\l|g^k\cdot \l[\kl\frac{g'}{g}\kr^{(k-1)}\r]^m\r|(z_0)\nonumber\\
&\le & \l|\frac{h}{B_*}(z_0)\r|^k\cdot 
\kl\frac{6\cdot k!\cdot 2^k}{(s-r)^{2k+1}} \cdot\log\frac{1}{|h(0)|}\kr^m
\cdot \kl m+\frac{1}{(s-r)^2}\cdot
\log\frac{|B_*(z_0)|}{|h(z_0)|}\kr^{2m}\nonumber\\
&\le & \l|h(z_0)\r|^k\cdot 
\kl\frac{6\cdot k!\cdot 2^k}{(s-r)^{2k+1}} \cdot\log\frac{1}{|h(0)|}\kr^m
\cdot \kl m+\frac{1}{(s-r)^2}\cdot \log\frac{1}{|h(z_0)|}\kr^{2m},
\eeqar
i.e. we have eliminated the function $B_*$ which depended on
$z_0$. We have shown this estimate for all $z_0\in\kj{\Delta(0,r)}$
with $g(z_0)\ne0$. By continuity, it even holds for all $z_0\in
\kj{\Delta(0,r)}$. (Note that the function $g^k\cdot
\l[\kl\frac{g'}{g}\kr^{(k-1)}\r]^m$ is holomorphic on $\dk$ by our
assumption on the multiplicities of the zeros of $g$.)

These considerations were subject to the assumption $n(s,0,g)\ge
1$. However, if $n(s,0,g)=0$, then $B\equiv 1$ and $g=h$, and from
\eqref{g-h-B} and \eqref{h-LogDeriv} we immediately obtain 
$$\l|g^k\cdot \l[\kl\frac{g'}{g}\kr^{(k-1)}\r]^m\r|(z)\\
\le  \l|h(z)\r|^k\cdot 
\kl\frac{2\cdot k!}{(s-r)^{k+1}} \cdot\log\frac{1}{|h(0)|}\kr^m \qquad
\mbox{ for } |z|\le r,$$
i.e.~an estimate even better than (\ref{B*Elim}). So in both cases
$n(s,0,g)\ge 1$ and $n(s,0,g)=0$, (\ref{B*Elim}) holds for all $z$
with $|z|\le r$. 

Now, by Harnack's inequality we have
$$\log\frac{1}{|h(0)|}\le \frac{s+r}{s-r}\cdot\log\frac{1}{|h(z)|} 
\qquad\mbox{ for } |z|\le r$$ 
and we finally conclude that for $|z|\le r$ 
\beqar\label{KeyEstim}
&&\l|g^k\cdot \l[\kl\frac{g'}{g}\kr^{(k-1)}\r]^m\r|(z)\nonumber\\
&\le & \l|h(z)\r|^k\cdot 
\kl\frac{12s\cdot k!\cdot 2^k}{(s-r)^{2k+2}} \cdot\log\frac{1}{|h(z)|}\kr^m
\cdot \kl m+\frac{1}{(s-r)^2}\cdot \log\frac{1}{|h(z)|}\kr^{2m}.
\eeqar
This estimate holds for all holomorphic functions $g$ on $\dk$ all of
whose zeros have multiplicity at least $m$  and which satisfy
$|g(z)|\le x_0$ for all $z\in \dk$.  

We apply this estimate to the sequence $\gl g_n\gr_n$. To each $g_n$,
as in (\ref{h-Def}) we construct a function $h_n$ holomorphic on
$\dk$ and non-vanishing in $\Delta(0,s)$ such that
$|h_n(z)|=|g_n(z)|$ for $|z|=s$, hence
$$\max_{|z|\le s} |h_n(z)| =\max_{|z|\le s} |g_n(z)| \abb0 \quad
(n\to\infty)$$ 
by the maximum principle. Now from (\ref{KeyEstim}) we immediately obtain
that the sequence $\gl g_n^k\cdot \l[\kl\frac{g_n'}{g_n}\kr^{(k-1)}\r]^m\gr_n$
converges to 0 uniformly in $\Delta(0,r)$. Since this holds for any
$r<1$, our assertion in I. is proved. 

{\bf II.} Now we prove (a) by induction. The  case $k=1$ follows immediately
from I.  

Let some $k\ge 2$ be given and assume that 
$$\frac{\kl g_n^{(j)}\kr^m}{g_n^{m-j}}\Longrightarrow 0 \qquad\mbox{
  for } j=1,\dots,k-1$$
has already been proved. Now by induction there are certain universal
  constants $c_{k;l;j_1,\dots,j_l}$ such that   
$$\frac{g^{(k)}}{g} 
=\l(\frac{g'}{g}\r)^{(k-1)} 
+ \sum_{l=2}^k \sum_{j_1+\dots+j_l=k \atop j_\mu\geq 1} c_{k;l;j_1,\dots,j_l} 
\cdot\prod_{\mu=1}^l \frac{g^{(j_\mu)}}{g}$$
for all functions $g\not\equiv0$ holomorphic on $\dk$. Setting
$$S_{n,k}:=\sum_{l=2}^k \sum_{j_1+\dots+j_l=k \atop j_\mu\geq 1} c_{k;l;j_1,\dots,j_l} 
\cdot\prod_{\mu=1}^l \frac{g_n^{(j_\mu)}}{g_n},$$
we obtain 
\beqar\label{induction1}
\l| g_n^k\cdot \kl\frac{g_n^{(k)}}{g_n}\kr^m\r|
&\le& \l| g_n^k \cdot \l[\kl\frac{g_n'}{g_n}\kr^{(k-1)}\r]^m\r|\\
&&+ \sum_{\sigma=0}^{m-1} {m\choose \sigma} \cdot
\l|\kl\frac{g_n'}{g_n}\kr^{(k-1)}\r|^\sigma \cdot |g_n|^{k\sigma/m}
\cdot |S_{n,k}|^{m-\sigma}\cdot |g_n|^{k(m-\sigma)/m}.\nonumber
\eeqar
Here, from I. we know
$$\l|\kl\frac{g_n'}{g_n}\kr^{(k-1)}\r|^\sigma \cdot |g_n|^{k\sigma/m}
\Longrightarrow 0 \quad (n\to\infty) \qquad \mbox{
for } \sigma=1,\dots,m-1,$$
and from 
$$|S_{n,k}|\cdot |g_n|^{k/m} 
\le \sum_{l=2}^k \sum_{j_1+\dots+j_l=k \atop j_\mu\geq 1} |c_{k;l;j_1,\dots,j_l}|
\prod_{\mu=1}^l \l|\frac{\kl g_n^{(j_\mu)}\kr^m}{g_n^{m-j_\mu}}\r|^{1/m}$$
and the induction hypothesis we deduce that 
$$|S_{n,k}|\cdot |g_n|^{k/m} \Longrightarrow 0 \quad (n\to\infty).$$
Inserting this into \eqref{induction1} and observing I. once more
yields  
$$g_n^k\cdot \kl\frac{g_n^{(k)}}{g_n}\kr^m\Longrightarrow 0 \quad
(n\to\infty),$$
as asserted. 

{\bf III.} We turn to the proof of (b). If we apply I. to the functions
$g_n:=1/f_n$ (all of whose zeros have multiplicity at least $p$), we
obtain that under the assumptions in (b) for all $k\ge 1$ we have 
$$\frac{1}{f_n^k}\cdot \l[\kl\frac{f_n'}{f_n}\kr^{(k-1)}\r]^p\Longrightarrow
0 \quad(n\to\infty).$$ 
From this we deduce (b) almost literally as in II. we have deduced (a)
from I. 
\ebew

Once Theorem \ref{ChenHuaGen} (b) is available, Theorem
\ref{LiXieGeneraliz} can be proved with the same method as in the
proof of Theorem \ref{LX} given in \cite{NevoPang}. For completeness, we
provide the details. 

\bbew{ of Theorem \ref{LiXieGeneraliz}}
We assume that $\fa$ is normal but that  $\fa_{k,\alpha}$ is not
locally uniformly bounded in $D$. Then we find a $z_0\in D$, functions
$f_n\in \fa$ and points $z_n\in D$ such that $\limn z_n=z_0$ and 
\beq\label{DerivUnbounded}
\dfrac{|f_n^{(k)}|}{1+|f_n|^\alpha}(z_n)\underset{n\to\infty}\abb\infty.
\eeq
Since $\fa$ is normal, after extracting a suitable subsequence we may
assume that $\gl f_n\gr_n$ converges locally uniformly to some limit
function $f$, possibly $f\equiv \infty$. Now let us consider several
cases. 

 \noindent \textbf{Case 1.}
$f(z_0)\in\mathbb C.$

Then there are $r>0$ and $N\in\nat$ such that $f$ and $f_n$ are holomorphic on
$\Delta(z_0,r)$ for all $n\ge N$, and by Weierstra\ss{}'s theorem we obtain
$$\limn\frac{f_n^{(k)}(z_n)}{1+|f_n(z_n)|^\alpha}=
\frac{f^{(k)}(z_0)}{1+|f(z_0)|^\alpha} \ne\infty,$$ 
a contradiction to (\ref{DerivUnbounded}).

\noindent \textbf{Case 2.} $f\not\equiv\infty$, but $f(z_0)=\infty.$

Here, we can find $r>0$ such that $f$ is holomorphic on
$\Delta'(z_0,2r)$ and $|f(z)|\ge 1$ and $|f_n(z)|\ge 1$ for all
$z\in\Delta(z_0,2r)$ and all $n\ge N$ for a certain $N\in\nat$.

If $p:=\lceil \frac{k}{\alpha-1}\rceil$ is the smallest integer
$\ge\frac{k}{\alpha-1}$, then, by assumption, each pole of $f_n$ has
multiplicity at least $p$.  Hence the functions 
\beq\label{Def-Dn}
d_n:=\frac{\kl f_n^{(k)}\kr^p}{f_n^{p+k}}
\eeq
are holomorphic on $\Delta(z_0,2r)$ for $n\ge N$. Since they converge to 
$\frac{\kl f^{(k)}\kr^p}{f^{p+k}}$ uniformly on
$\partial\Delta(z_0,r)$, from the maximum principle we deduce that
there is a constant $C<\infty$ such that 
$$|d_n(z)|\le C \qquad\mbox{ for all } z\in \Delta(z_0,r) \mbox{ and $n$
  large enough}.$$
In particular, for $z=z_n$ we get for $n$ large enough 
$$\kl\frac{|f_n^{(k)}(z_n)|}{1+|f_n(z_n)|^\alpha}\kr^p
\le\frac{|f_n^{(k)}(z_n)|^p}{|f_n(z_n)|^{p+k}} \le C;$$
here we have used $|f_n(z_n)|\ge 1$ and $\alpha p \ge
k+p$. This is a contradiction to (\ref{DerivUnbounded}).  

\noindent \textbf{Case 3.} $f  \equiv\infty.$

Again, each pole of $f_n$ has multiplicity at least $p:=\lceil
\frac{k}{\alpha-1}\rceil$. So from Theorem \ref{ChenHuaGen} (b) we
obtain that the sequence $\gl d_n\gr_n$ where $d_n$ is defined as in
(\ref{Def-Dn}) converges to 0 locally uniformly in $D$, and in view of
$k+(1-\alpha)p\le k+(1-\alpha)\cdot\frac{k}{\alpha-1}=0$ we deduce
$$\dfrac{|f_n^{(k)}|}{1+|f_n|^\alpha}
\le \kl |d_n|\cdot |f_n|^{k+(1-\alpha)p}\kr^{1/p}\Longrightarrow 0$$
for $n\to\infty$, a contradiction to (\ref{DerivUnbounded}) once more. 
\ebew

\bibliographystyle{amsplain}

\end{document}